\documentclass[12pt]{amsart}
\numberwithin{equation}{section}
\newtheorem{theorem}{Theorem}
\newtheorem{lemma}{Lemma}

\numberwithin{theorem}{section} \numberwithin{lemma}{section}

\numberwithin{proposition}{section}
\begin{document}

\tracingpages 1
\title[blow up]{ finite time blow up for critical
wave equations in high dimensions}
\author{Borislav T. Yordanov and Qi S. Zhang}
\address{Department of Mathematics\\
University of California Riverside, Riverside, CA 92521}
\date{March 2004}

\begin{abstract}
We prove that solutions to the critical wave equation (1.1) can
not be global if the initial values are positive somewhere and
nonnegative. This completes the solution to the famous blow up
conjecture about critical semilinear wave equations of the form
$\Delta u - \partial^2_t u + |u|^p = 0$ in dimensions $n \ge 4$.
The lower dimensional case $n \le 3$ was settled many years
earlier.
\end{abstract}
\maketitle

\section{ Introduction}
Let $n\geq 2$ and $\Delta=\Sigma_{i=1}^{n}\partial^2/\partial
x_i^2$ be the Laplace operator. We study the blow up of solutions
to the following semilinear wave equation:
\[
\left\{
\begin{array}{ll}\Delta u - \partial^2_t u + |u|^p=0\quad
                 {\rm in}\quad {\bf R}^n\times(0,\infty),\\
                 u(x,0)=u_0(x), \quad u_t(x,0)=u_1(x)\quad
                 {\rm in}\quad {\bf R}^n,
\end{array}
\right. \leqno(1.1)
\]
where the initial values satisfy
\[
\left\{
\begin{array}{ll}
 (u_0,u_1)\in H^1({\bf R}^n)\times L^2({\bf
                 R}^n),\\
                  u_0(x)=u_1(x)=0 \quad {\rm for} \quad |x|>R>0,
\end{array}
\right.
\]
and $p\in(1,\ p_c(n)].$ Here $p_c(n)$ is the positive root of the
quadratic equation
\[
(n-1)p^2-(n+1)p-2=0.
\]
The number $p_c(n)$ is known as the critical exponent of problem
(1.1), since it divides $(1,\ \infty)$ into two subintervals so
that the following take place: If $p\in(1,\ p_c(n)),$ then
solutions with nonnegative initial values blow up in finite time;
if $p\in(p_c(n), \ \infty),$ then solutions with small (and
sufficiently regular) initial values exist for all time (see [St]
e.g.). The proof has an interesting and exciting history that
spans three decades. We only give a brief summary here and refer
the reader to [St], [L], [DL] and a recent paper [JZ] for details.
The problem about existence or nonexistence of global solutions is
sometimes referred to as the conjecture of Strauss [St2]. The
question was also asked by Glassey [G2] and Levine [L].

The case $n=3$ was considered by John [J] who proved that
nontrivial solutions must blow up in finite time when
$1<p<p_c(3)$. He also showed that global solutions exist for small
initial values when $p>p_c(3)$. Glassey [G1], [G2] established the
same results in the case $n=2$. In [GLS] Georgiev, Lindblad and
Sogge showed the existence of global solutions for small initial
values when $p>p_c(n)$ and $n \ge 4$. (See also the work of Kubo
and Kubota [KK],  Lindblad and Sogge [LS] and Tataru[T]). The
corresponding blow up result for $1<p<p_c(n)$ and $n \ge 4$ was
established by Sideris [Si]. Although the ideas of [Si] are very
clear, the computations are quite sophisticated, involving
spherical harmonics and other special functions. The proof was
simplified by Ramaha [R] and Jiao and Zhou [JZ].

The critical case $p=p_c(n)$ was studied by Schaeffer [Sc] in
dimensions $n=2$ and $3$. Improving the lower bounds on the
solution in [G2] and [Si], he was able to show that the critical
powers belong to the corresponding blow up intervals. Despite the
long effort, whether the critical powers $p_c(n)$ belong to the
blow up intervals remains wide open in dimensions $n \ge 4$. The
main obstruction to the method of [Sc] is that the Riemann
function changes sign in high dimensions. This difficulty is not
present if the initial values are large; the work of Levine [L]
shows that such solutions blow up in finite time. Thus, the open
problem is to show blow up without the latter assumption.

Here, we complete the solution of this conjecture about equation
(1.1), thus filling the missing link since the 80s. Our main
result is the following theorem.

\begin{theorem}
Let $u_0$ and $u_1$ be non-negative and let either of them be
positive somewhere. Suppose that problem (1.1) has a solution $(u,
u_t)\in C([0,T),\ H^1({\bf R}^n)\times L^2({\bf R}^n))$ such that
\[
supp(u,u_t)\subset \{(x,t):|x|\leq t+R \}.
\]
If $p=p_c(n)$, then $T<\infty$. \label{th1}
\end{theorem}
\medskip

We should mention that the existence of local in time solutions
with the above regularity and support is well known. See p381 in
[Si], for example. We prove Theorem 1.1 in Section~2. The cases
$n=2$ and $3$ are proven in [Sc], so we concentrate on the case $n
\ge 4$. Following the tradition, we consider $\int_{{\bf R}^n}
u(x, t) dx,$ where $u$ is a local solution of problem (1.1). We
show that this quantity satisfies a nonlinear differential
inequality and, additionally, admits a lower bound
$O(K(t)t^{n+1-(n-1)p/2})$ with $K(t)\geq \ln t$ as
$t\rightarrow\infty$. The finite time blow up then follows
immediately. Our estimate improves $K(t)\geq 1$, which is
sufficient to show blow up only in the subcritical case. The new
tools used are the Radon transform and the one-dimensional
transform ${\bf T}$ (see (2.16)). These together with the $L^p$
boundedness of the maximal function yield the extra $\ln t$ factor
in our lower bound.

\section{Proof of Theorem~\ref{th1}}

The proof is carried out in several steps. We assume that
$p=p_c(n)$ and $n \ge 4$.
\medskip

{\bf Step 1.}
\medskip

 We will need the
following ODE result.

\begin{lemma}
Let $p>1$, $a \ge 1$, and $(p-1) a=q-2$. Suppose $F\in C^2([0,T))$
satisfies, when $t \ge T_0>0$,
\[
\aligned
(a)& \qquad F(t)  \ge  K_0 (t+R)^a, \\
(b)& \qquad \frac{d^2 F(t)}{dt^2}  \ge  K_1(t+R)^{-q}[F(t)]^p,
\endaligned
\]
with some positive constants $K_0$, $K_1$, $T_0$ and $R$. Fixing
$K_1$, there exists a positive constant $c_0$, independent of $R$
and $T_0$ such that if $K_0 \ge c_0$, then $T<\infty$. \proof
\end{lemma}

First let us make a translation $\tau = t-T_0$ and define
$G=G(\tau)=F(\tau+R)$. Then for $\tau \ge 0$, one has
\[
 G(\tau)  \ge  K_0 (\tau +T_0 +R)^a, \qquad
 \frac{d^2 G(\tau)}{d\tau^2}  \ge
K_1(\tau+T_0+R)^{-q}[G(\tau)]^p.
\]

We take the change of variables $\tau = (T_0+R)  s$ and $G_R =
G_R(s)= (T_0+R)^{-a} G((T_0+R) s)$. Then easy computation shows
that
\[
G_R(s)  \ge  K_0 (s+1)^a, \qquad \frac{d^2 G_R(s)}{ds^2}  \ge
K_1(s+1)^{-q}[G_R(s)]^p,
\]when $s \ge 0$. Following the argument in [Si], p386, we know
that $G_R$ has to blow up in finite time if $K_0 \ge c_0$, which
is sufficiently large. Clearly $c_0$ does not depend on $R$ or
$T_0$. Therefore $F$ must also blow up in finite time. \qed

\medskip

{\bf Step 2.}
\medskip

We introduce the function
\[
\phi_1(x)=\int_{S^{n-1}}e^{x\cdot\omega} d\omega. \]It is well
known that
\[
\phi_1(x)\sim C_n |x|^{-(n-1)/2}{e^{|x|}}\quad {\rm
as} \quad |x|\rightarrow \infty.
\]

Suppose (1.1) has a global solution under the given initial
values. Define
\[
\aligned
F_0(t)&=\int u(x,t) dx,\\
F_1(t)&=\int u(x,t)\psi_1(x,t)dx, \\
 \psi_1(x&,t)=\phi_1(x)e^{-t}.
\endaligned
\leqno(2.2)
\]

To show that $F_0$ satisfies the  differential inequality in Lemma
2.1 for suitable $a$, $q$, we integrate equation (1.1) over ${\bf
R}^n$. We know that the support of $u(\cdot, t)$ is contained in
$B(0, t+R)$ since the supports of $u_0, u_1$ are contained in
$B(0, R)$. Hence using integration by parts, we obtain

\[
\frac{d^2 F_0(t)}{dt^2}=\int |u(x,t)|^p dx. \leqno(2.2')
\]

Estimating the right side by the H\"{o}lder inequality, we have
\[
\int |u(x, t)|^p dx \geq \frac{\left|\int u(x, t) dx
\right|^p}{\left(\int_{|x|\leq t+R} dx \right)^{p-1}}.
\]
Since
\[
\int_{|y|\leq t+R} dx =  {\rm vol}\{x: |x|<t+R \} = {\rm vol}({\bf
B}^n)(t+R)^n,
\]
we obtain the differential inequality
\[
\frac{d^2 F_0(t)}{dt^2} \geq K_1(t+R)^{-n(p-1)}|F_0(t)|^p
\leqno(2.3)
\]
with  $K_1=1/({\rm vol} ({\bf B}^n))^{p-1}.$

To show that $F_0$ admits the lower bound in Lemma 2.1 (a), we
relate $d^2F_0/dt^2$ to $F_1$ using again equation (1.1) and
H\"{o}lder's inequality:
\[
\frac{d^2F_0(t)}{dt^2}=\int |u(x, t)|^p dx \geq \frac{\left|\int
u(x, t)\psi_1(x, t)dy \right|^p}{\left(\int_{|x|\leq t+R}
[\psi_1(x, t)]^{p/(p-1)}dx\right)^{p-1}}.
\]By (2.2), the above becomes
\[
\frac{d^2F_0(t)}{dt^2} \ge \frac{ \left|
F_1(t)\right|^p}{\left(\int_{|x|\leq t+R}
[\psi_1(x,t)]^{p/(p-1)}dx\right)^{p-1}}. \leqno(2.4)
\]
In the  following we estimate the denominator and numerator,
respectively.

We claim that for all $ t \geq 0$, $R>0$,
\[
I(t) \equiv \int_{|x| \leq t+R} [\psi_1(x,t)]^{p/(p-1)} dx \leq C
e^{p' R} (t+R)^{n-1-(n-1)p'/2}, \leqno(2.5)
\]
where $p'=p/(p-1)$. The claim is an immediate consequence of the
observation
\[
I(t) \le C_1 e^{-p't} \int_0^{t+R}(1+r)^{-(n-1)p'/2}
e^{p'r}r^{n-1} dr,
\]with $p'=p/(p-1)$ and integration by parts. Here we just used
the formula
\[
\psi_1(x, t) = e^{-t} \phi_1(x)\sim C_n
|x|^{-(n-1)/2}{e^{|x|-t}}\quad {\rm as} \quad |x|\rightarrow
\infty.
\]

\medskip

Next we have

\begin{lemma}
For all $t \geq 0$,
\[
F_1(t) \ge  \frac{1}{2}(1-e^{-2t})\int [u_0(x)+ u_1(x)] \phi_1(x)
dx
      + e^{-2t} \int u_0(x)\phi_1(x)dx \ge c>0.
\]
\end{lemma}
\medskip

Taking the lemma for granted, we combine it with (2.5) and with
(2.4) to obtain
\[
\frac{d^2F_0(t)}{dt^2} = \int_{{\bf R}^n} |u(x, t)|^p dx \geq C_0
L_2 (t+R)^{n-1-(n-1)p/2}, \quad t\geq 0, \leqno(2.5')
\]
where
\[
L_2 \ge \bigg{(} C \int u_0(x) \phi_1(x) dx \bigg{)}^p, \qquad
C>0. \leqno(2.6)
\]
Integrating twice, we have the  estimate
\[
F_0(t)\geq c L_2 (t+R)^{n+1-(n-1)p/2} +\frac{dF_0(0)}{dt} t+F_0(0)
\]
with some $c>0$ depending only on $n$. When $p=p_c(n)$, it is easy
to check that $n+1-(n-1)p/2>1$. Hence the following estimate is
valid for all $t \ge 0$:
\[
F_0(t) \ge K_0 (t+R)^{n+1-(n-1)p/2}. \leqno(2.7)
\]with $K_0 \equiv c L_2$.
Here we remark that (2.7) have been proven in [Si] and [JZ] by
different method. The current method, adopted from [YZ], seems
much shorter.

If $K_0$ is sufficiently large, estimates (2.7), (2.3), and Lemma
2.1 with parameters
\[
a \equiv n+1-(n-1)p/2 \quad {\rm and} \quad  q \equiv n(p-1)
\]would
imply Theorem 1.1 since $p=p_c$ satisfies
\[
(p-1)(n+1-(n-1)p/2)=n(p-1)-2 \quad {\rm and} \quad p>1.
\]However we have no control on the size of $K_0$.  In the
remainder of the paper, we will show that the lower bound (2.7)
can be improved by a factor of $\ln t$ when $t$ is large. Before
doing so let us give a

\medskip

{\it Proof Lemma 2.2.}

  We multiply equation (1.1) by a test
function $\psi\in C^2({\bf R}^{n+1})$ and integrate over ${\bf
R}^n\times [0,\, t].$
\[
\aligned
  \int_0^t\int & u \ (\Delta \psi
  - \partial^2_s \psi)dyds
      +\int_0^t \int |u|^p \psi \, dyds\\
& =  \int(\partial_s u \ \psi - u
\partial_s \psi ) dy \big{|}^{s=t}_{s=0}.
\endaligned
\leqno(2.8)
\]
We will apply this identity to $\psi=\psi_1$. Notice that for a
fixed $t$, $u(\cdot, t) \in H^1_0(B(0, t+R))$. Hence all terms
involving lateral boundary vanish during integration by parts.
Notice also that
\[
\partial_t \psi_1=-\psi_1, \quad \Delta\psi_1 -\partial^2_t
\psi_1 =0,
\]
and
\[
\aligned \int(\partial_t u \psi_1-u
\partial_t \psi_1) dy &=\int( \partial_t u \psi_1
+u \partial_t \psi_1) dy
-2\int u \partial_t \psi_1 dy \\
& =\frac{d}{dt}\int u \psi_1 dy +2\int u \psi_1 dy.
\endaligned
\]

Hence, (2.8) becomes
\[
\frac{dF_1(t)}{dt}+2F_1(t)=\int[u(x, 0)+
\partial_t u(y, 0)]\phi_1(y)dy+ \int_0^t\int|u(y,s)|^p
\psi_1(y,s) dyds.
\]
Since $\psi_1>0,$ we conclude that
\[
\frac{dF_1(t)}{dt}+2F_1(t) \geq \int [u(y, 0)+
\partial_t u(y, 0)]\phi_1(y)dy.
\]
We multiply by $e^{2t}$ and integrate on $[0,\, t]$. Then
\[
e^{2t}F_1(t)- F_1(0)\geq \frac{1}{2}(e^{2t}-1) \int[u_0(y)+
u_1(y)]\phi_1(y)dy.
\]
Dividing through by $e^{2t}$, we obtain the lower bound in the
Lemma. \qed

\medskip

{\bf Step 3.}
\medskip

With no loss of generality we assume that $u(\cdot, t)$ is radial.
This is so because one can use Daboux's identity to transform the
problem into a suitable inequality in the radial case. i.e the
sperical average of $u$, called $\bar{u}$ satisfies
\[
\partial^2_t \bar{u} - \Delta \bar{u} \ge |\bar{u}|^p.
\]

Let $w \in {\bf R}^n$ be a unit vector. The Radon transform of $u$
with respect to the space variables is defined as
\[
{\bf R}(u)(\rho, t) = \int_{x \cdot w = \rho} u(x, t) dS_x,
\leqno(2.9)
\]where $dS_x$ is the Lebesque measure on the hyper-plane
$\{ x \ | \ x \cdot w = 0 \}.$  Next we show that ${\bf R}(u)$ is
 a function of $\rho$ and $t$ and is in fact independent of $w$.

From (2.9) and the assumption that $u(\cdot, t)$ is radial, it is
clear that
\[
\aligned
 {\bf R}(u)(\rho, t) &= \int_{ \{ x' \ | \ x' \cdot w = 0
\} } u(\rho w + x', t) dS_{x'} \\
&= c_n \int^{\infty}_0 u(\sqrt{\rho^2+|x'|^2}, t) |x'|^{n-2}
d|x'|.
\endaligned
\]Using the change of variable $r = \sqrt{\rho^2+|x'|^2}$, we have
\[
{\bf R}(u)(\rho, t)  = c_n \int^{\infty}_{|\rho|} u(r, t)
(r^2-\rho^2)^{(n-3)/2} \ r dr. \leqno(2.10)
\]This shows that ${\bf R}(u)(\rho, t)$ is independent of $w$.
In the remainder of the step, we will derive a lower  bound for
${\bf R}(u)(\rho, t)$.

Since $u$ is a solution to (1.1), it is well known that ${\bf
R}(u)$ satisfies the one dimensional wave equation
\[
\partial^2_t {\bf R}(u)(\rho, t) -
\partial^2_{\rho} {\bf R}(u)(\rho, t) =
{\bf R}( |u|^p )(\rho, t). \leqno(2.11)
\]From the D' Alembert's formula and the assumption that the
initial values of $u$ are nonnegative, one obtains
\[
{\bf R}(u)(\rho, t) \ge \frac{1}{2} \int^t_0 \int^{\rho +
(t-s)}_{\rho - (t-s)} {\bf R}( |u|^p )(\rho_1, s) d\rho_1ds.
\leqno(2.12)
\]

Observe that the support of $u(\cdot, s)$ is contained in $B(0,
s+R)$, the ball of radius $R$, centered at the origin. If
$|\rho_1|> s+R$, then, for vectors $y$ perpendicular to a unit
vector $w$,
\[
|\rho_1 w + y| = \sqrt{|\rho_1|^2 +|y|^2} \ge |\rho_1|> s+R.
\]Therefore
\[
{\bf R}( |u|^p )(\rho_1, s) = \int_{\{y \ | y \cdot w = 0\}}
|u(\rho_1 w + y, s)|^p dS_y =0.
\]This shows that
\[
supp \ {\bf R}(|u|^p)(\cdot, s) \subset B(0, s+R). \leqno(2.13)
\]

From now on we will assume $\rho \ge 0$, unless stated otherwise.
If $s \le (t-\rho-R)/2$, then
\[
\rho + (t-s) \ge s+R, \qquad \rho-(t-s) \le -(s+R).
\]Using this, (2.12) and (2.13), we deduce
\[
\aligned
 {\bf R}(u)(\rho, t) &\ge \frac{1}{2} \int^{(t-\rho-R)/2}_0 \int^{\rho + (t-s)}_{\rho -
(t-s)} {\bf R}( |u|^p )(\rho_1, s)
d\rho_1ds\\
&= \frac{1}{2} \int^{(t-\rho-R)/2}_0 \int^{\infty}_{-\infty} {\bf
R}(
|u|^p )(\rho_1, s) d\rho_1ds\\
&=\frac{1}{2} \int^{(t-\rho-R)/2}_0 \int_{{\bf R}^n} |u(y, s)|^p
dyds.
\endaligned
\leqno(2.14)
\]Recall from (2.5$'$) in step 2 that
\[
\int_{{\bf R}^n} |u(y, s)|^p dy \ge c s^{(n-1) - (n-1)p/2}.
\]Note that $p \le 2$ when $n \ge 4$. Therefore $(n-1) - (n-1)p/2
\ge 0$.

Substituting this to (2.14), we arrive that
\[
{\bf R}(u)(\rho, t) \ge c (t-\rho - R)^{n - (n-1)p/2}, \qquad \rho
\ge 0. \leqno(2.15)
\]
\medskip

{\bf Step 4.} For any function $f\in L^p({\bf R})$, we introduce
the transformation
\[
{\bf T}(f)(\rho) = \frac{1}{|t-\rho+R|^{(n-1)/2}}
\int^{t+R}_{\rho} f(r) |r - \rho|^{(n-3)/2} dr. \leqno(2.16)
\]Observe that
\[
\aligned
 |{\bf T}(f)(\rho)|& \le \frac{1}{|t-\rho+R|}
\bigg{|} \int^{t+R}_{\rho} |f(r)| dr \bigg{|} \\
&\le \frac{2}{2|t-\rho+R|} \bigg{|} \int^{t+R}_{-(t+R)+2\rho}
|f(r)| dr \bigg{|} \\
&\le 2 M(|f|)(\rho),
\endaligned
\]where $M(|f|)$ is the maximal function of $f$.
Therefore, there exists a $C>0$, independent of $t$ such that
\[
\Vert {\bf T}(f) \Vert_p \le C \Vert f \Vert_p. \leqno(2.17)
\]Here we remark that (2.17) can also be proven directly by
showing that ${\bf T}$ maps $L^\infty$ to $L^{\infty}$ and $L^1$
to weak $L^1$. Then the Marcinkiewicz interpolation theorem will
imply (2.17).

Applying (2.17) to the function \[ f(r)= \begin{cases} |u(r, t)|
r^{(n-1)/p}, \qquad r \ge 0\\
0, \qquad r<0, \end{cases}
\]we have
\[
\aligned & \int^{t+R}_0 \bigg{[} \frac{1}{(t-\rho+R)^{(n-1)/2}}
\int^{t+R}_{\rho} |u(r, t)| r^{(n-1)/p} (r - \rho)^{(n-3)/2} dr
\bigg{]}^p d\rho \\
& \le C \int^\infty_0 |u(r, t)|^p r^{n-1} dr\\
& = C \int_{{\bf R}^n} |u(x, t)|^p dx.
\endaligned
\leqno(2.18)
\]When $r \ge \rho$ and $1<p \le 2$, we observe that
\[
r^{(n-1)/p} = r^{(n-1)/2} \ r^{(n-1)/p - (n-1)/2} \ge r^{(n-1)/2}
\ \rho^{(n-1)/p - (n-1)/2}.
\]Hence (2.18) becomes
\[
\aligned & \int^{t+R}_0 \bigg{[} \frac{1}{(t-\rho+R)^{(n-1)/2}}
\int^{t+R}_{\rho} |u(r, t)| r^{(n-1)/2} (r - \rho)^{(n-3)/2} dr
\bigg{]}^p \ \rho^{n-1 - (n-1)p/2} d\rho \\
& \le  C \int_{{\bf R}^n} |u(x, t)|^p dx.
\endaligned
\leqno(2.19)
\]

From (2.10) and the fact that $supp \  u(\cdot, t) \subset B(0,
t+R)$, we know that
\[
\aligned
 {\bf R}(|u|)(\rho, t) & = c_n \int^{t+R}_\rho |u(r, t)| r
(r^2-\rho^2)^{(n-3)/2} dr\\
&\le c_n \int^{t+R}_\rho |u(r, t)| r (r+\rho)^{(n-3)/2}
(r-\rho)^{(n-3)/2} dr\\
&\le c \int^{t+R}_\rho |u(r, t)| r^{(n-1)/2} (r-\rho)^{(n-3)/2}
dr.
\endaligned
\leqno(2.20)
\]Substituting (2.20) to (2.19), we reach
\[
\int^{t+R}_0  \frac{[{\bf R}(|u|)(\rho,
t)]^p}{(t-\rho+R)^{(n-1)p/2}}
 \ \rho^{n-1 - (n-1)p/2} d\rho \\
 \le  C \int_{{\bf R}^n} |u(x, t)|^p dx.
\leqno(2.21)
\]

Using the lower bound of ${\bf R}(|u|)$ in (2.15) and (2.21), we
deduce
\[
\int_{{\bf R}^n} |u(x, t)|^p dx \ge C \int^{t-R-1}_0 \frac{(t-\rho
- R)^{np - (n-1)p^2/2}}{(t-\rho+R)^{(n-1)p/2}}
 \ \rho^{n-1 - (n-1)p/2} d\rho.
\]When $\rho \in (0, t-R-1)$, it is clear that there exists
$c_R>0$ such that, for all $t>2(R+1)$,
\[
t-\rho+R \le c_R (t - \rho -R).
\]Hence there exist $C_R>0$ such that
\[
\int_{{\bf R}^n} |u(x, t)|^p dx \ge C_R \int^{t-R-1}_0
\frac{\rho^{n-1 - (n-1)p/2}}{(t-\rho-R)^{(n-1)p/2 - np +
(n-1)p^2/2}} d\rho. \leqno(2.22)
\]Recall that $p$ is the critical exponent for (1.1), i.e.
\[
(n-1) p^2 - (n+1) p - 2 = 0.
\]It follows that
\[
(n-1)p/2 - np + (n-1)p^2/2 = \frac{(n-1) p^2 - (n+1) p}{2}=1.
\]Therefore (2.22) becomes
\[
\int_{{\bf R}^n} |u(x, t)|^p dx \ge C_R \int^{t-R-1}_0
\frac{\rho^{n-1 - (n-1)p/2}}{(t-\rho-R)} d\rho. \leqno(2.23)
\]Hence
\[
\aligned \int_{{\bf R}^n} &|u(x, t)|^p dx \\
&\ge C_R \int^{t-R-1}_{(t-R-1)/2}
\frac{\rho^{n-1 - (n-1)p/2}}{(t-\rho-R)} d\rho\\
&\ge  C_R (t-R-1)^{n-2 - (n-1)p/2} \int^{t-R-1}_{(t-R-1)/2}
\frac{1}{t-\rho-R} d\rho.
\endaligned
\]So finally, we reach the refined lower bound
\[
\int_{{\bf R}^n} |u(x, t)|^p dx \ge C (t-R)^{n-1 - (n-1)p/2}
\ln\frac{t-R}{2}. \leqno(2.24)
\]From (2.2'), the above shows
\[
\frac{d^2 F_0(t)}{dt^2} = \int_{{\bf R}^n} |u(x, t)|^p dx \ge C
(t-R)^{n-1 - (n-1)p/2} \ln\frac{t-R}{2},
\]provided that $t$ is sufficiently large. Comparing with the
lower bound in (2.5$'$), the above contains an additional $\ln t$
term. This the key improvement.

Since $n-1 - (n-1)p/2 \ge 0$ when $n \ge 4$, after integration we
deduce, for large $t$,
\[
F_0(t) \ge C (t-R)^{n+1 - (n-1)p/2} \ln t.
\]Hence
\[
F_0(t) \ge C (t+R)^{n+1 - (n-1)p/2} \big{(}
\frac{t-R}{t+R}\big{)}^{n+1 - (n-1)p/2}\ln t,
\]when $t$ is sufficiently large. Notice that
\[
\lim_{t \to \infty} \big{(} \frac{t-R}{t+R}\big{)}^{n+1 -
(n-1)p/2}\ln t = \infty.
\]Therefore
\[
F_0(t) \ge K_0 (t+R)^{n+1 - (n-1)p/2} \leqno(2.25)
\]with $K_0>0$ being arbitrarily large when $t$ is sufficiently
large.

Also, recall from (2.3) that
\[
\frac{d^2 F_0(t)}{dt^2} \geq K_1(t+R)^{-n(p-1)}|F_0(t)|^p
\]with $K_1 =1/({\rm vol} ({\bf B}^n))^{p-1}$ being fixed.

This together with (2.25) and Lemma 2.1 with parameters
\[
a \equiv n+1-(n-1)p/2 \quad {\rm and} \quad  q \equiv n(p-1)
\]
imply Theorem 1.1 since $p=p_c$ satisfies
\[
(p-1)(n+1-(n-1)p/2)=n(p-1)-2 \quad {\rm and} \quad p>1.
\]
This shows that all solutions of (1.1) with nontrivial nonnegative
initial values must blow up in finite time. \qed

\bigskip

{\bf Acknowledgement.} We thank Professors Ken Deng, Howard A.
Levine and Thomas Sideris for helpful conversations.

\bigskip

\noindent e-mail: yordanov@math.ucr.edu and  qizhang@math.ucr.edu

\enddocument